\documentclass[preprint, 3p, table]{elsarticle}

\makeatletter
\def\ps@pprintTitle{%
		\let\@oddhead\@empty
		\let\@evenhead\@empty
		\def\@oddfoot{}%
		\let\@evenfoot\@oddfoot}
\makeatother

\usepackage[english]{babel}
\usepackage[dvipsnames]{xcolor}
\usepackage[colorlinks=true,linkcolor=black,anchorcolor=black,citecolor=black,filecolor=black,menucolor=black,runcolor=black,urlcolor=black]{hyperref} 
\usepackage{changepage} 
\usepackage{amsmath}
\usepackage{amssymb}
\usepackage{lineno} 

\usepackage{graphicx}

\usepackage{amsthm} 
\usepackage{mathtools}
\usepackage{mathrsfs}
\usepackage{upgreek}
\usepackage{siunitx}
\usepackage{svg}

\usepackage{gensymb}
\usepackage{comment}
\usepackage[colorinlistoftodos]{todonotes}

\usepackage{algorithm}
\usepackage{algorithm,algpseudocode}
\usepackage{physics}
\usepackage{acro}

\usepackage{subfig}

\usepackage{tikz}

\usepackage{tabularx}
\usepackage{multirow}
\usepackage{booktabs}

\usepackage{empheq}

\everymath{\displaystyle}

\usepackage{hhline}

\usepackage[toc,page]{appendix}
\usepackage[capitalise,nameinlink]{cleveref}
\crefname{supp}{Supplement}{Supplements}

\usepackage{wasysym} 
\newif\ifdouble

\newcommand{\papertitle}{
Coupled Eikonal problems to model cardiac reentries \\ in Purkinje network and myocardium
}

\newcommand{\keywordOne}{Cardiac electrophysiology}
\newcommand{\keywordTwo}{Purkinje network}
\newcommand{\keywordThree}{Electrical reentries}
\newcommand{\keywordFour}{Eikonal models}
\newcommand{\keywordFive}{Pseudo-time method} 
\newcommand{\keywordSix}{Bundle branch block}

\definecolor{lifex}{HTML}{f60248}
\newcommand{\lifex}{\texttt{life\textsuperscript{\color{lifex}{x}}}}
\newcommand{\dutchcal}[1]{\text{\large\usefont{U}{dutchcal}{m}{n}#1}}

\DeclareAcronym{PMJs}{
	short = PMJs,
	long = Purkinje-muscle junctions
}
\DeclareAcronym{CCS}{
	short = CCS,
	long = cardiac conduction system
}
\DeclareAcronym{CRT}{
	short = CRT,
	long = cardiac resynchronization therapy
}
\DeclareAcronym{ECGs}{
	short = ECGs,
	long = electrocardiograms
}
\DeclareAcronym{FMM}{
	short = FMM,
	long = Fast Marching method
}
\DeclareAcronym{AV}{
	short = AV,
	long = atrioventricular
}
\DeclareAcronym{BDF}{
	short = BDF,
	long = backward differentiation formula
}
\DeclareAcronym{CT}{
	short = CT,
	long = computed tomography
}

\DeclareMathAlphabet{\mathdutchcal}{U}{dutchcal}{m}{n}
\graphicspath{{}}
\doublefalse

\begin{document}	
	
	\begin{frontmatter}
		\title{{\papertitle}}
		
		\journal{Journal of Computational Physics}

\author[mox]{Samuele~Brunati}
\author[mox]{Michele~Bucelli}
\author[mox]{Roberto~Piersanti\corref{cor1}}\ead{roberto.piersanti@polimi.it}
\author[mox]{Luca~Dede'}
\author[labs]{Christian Vergara}

\address[mox]{MOX, Laboratory of Modeling and Scientific Computing, Dipartimento di Matematica, Politecnico di Milano, Milano, Italy}
\address[labs]{LaBS, Laboratory of Biological Structure Mechanics, Dipartimento di Chimica, Materiali e Ingegneria Chimica, Politecnico di Milano, Milano, Italy}

\cortext[cor1]{Corresponding author.}
		
		\begin{abstract}
			We propose a novel partitioned scheme based on Eikonal equations to model the coupled propagation of the electrical signal in the His-Purkinje system and in the myocardium for cardiac electrophysiology. This scheme allows, for the first time in Eikonal-based modeling, to capture all possible signal reentries between the Purkinje network and the cardiac muscle that may occur under pathological conditions. As part of the proposed scheme, we introduce a new pseudo-time method for the Eikonal-diffusion problem in the myocardium, to correctly enforce electrical stimuli coming from the Purkinje network. We test our approach by performing numerical simulations of cardiac electrophysiology in a real biventricular geometry, under both pathological and therapeutic conditions, to demonstrate its flexibility, robustness, and accuracy.
		\end{abstract}
		
		
	\end{frontmatter}
	
	\section{Introduction}
	\label{sec:introduction}
\begin{figure}
	\centering
	\includegraphics[width=1\textwidth]{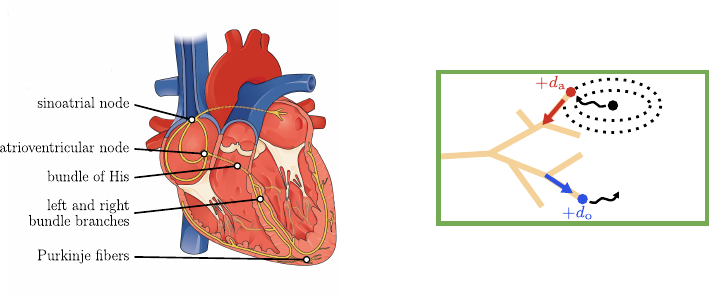}
	\caption{Left: representation of the components of the \acl{CCS} (picture adapted from \protect\url{https://commons.wikimedia.org/wiki/File:2018_Conduction_System_of_Heart.jpg}). Focusing the attention on the ventricular electrical activity, the signal coming from the atrioventricular node (AV node) descends into the bundle of His, which divides into the right and the left bundle branches, finally branching into the Purkinje network. Right: a visualization of orthodromic (in blue) and antidromic propagation (in red) in the Purkinje network (represented in yellow) due to an intramuscular source (in black). The signal undergoes an orthodromic delay $d_\mathrm{o}$ before being transmitted from the Purkinje to the muscle, and an antidromic delay~$d_\mathrm{a}~<~d_\mathrm{o}$ in the opposite direction.}
	\label{fig:intro}
\end{figure}

	Heart contraction is triggered by an electrical action potential propagating throughout the cardiac muscle. This process is regulated by the \ac{CCS}, schematized in Figure~\ref{fig:intro}(left). The Purkinje network plays a pivotal role in driving the action potential~\cite{iaizzo,quarteroni103,quarteroni2019mathematical,fisiologia_ita}. It is a dense network of conductive fibers, lying on the endocardium, and electrically connected to the surrounding tissue at its terminal points, called \textit{Purkinje-muscle junctions} (PMJs). Propagation is referred to as \textit{orthodromic} when excitation travels from the Purkinje network to the muscle, or \textit{antidromic} when the opposite happens, as depicted in Figure~\ref{fig:intro}(right). Antidromic propagation is typically absent in healthy hearts; however, it becomes relevant under diseased conditions~\cite{romano2010,palamara_anti,vergara2016}, where \textit{reentries}\footnote{Cardiological literature generally refers to reentry as an abnormal pathway for heart signals to travel, forming a reentrant circuit, regardless of whether it involves the His-Purkinje system~\cite{fleischmann1976reentry}. In this work, we focus exclusively on reentrant circuits involving the Purkinje network. With a slight misuse of terminology, we refer to such circuits simply as reentries.} in the network may occur. Therefore, incorporating both orthodromic and antidromic propagation mechanisms in computational models of cardiac electrophysiology is crucial for accurately reproducing the physical, possibly pathological, processes.

Existing cardiac computational models often lack an explicit representation of the Purkinje network, and usually surrogate it through realistic stimulation protocols based on sparse endocardial sources~\cite{vigmond2003,nickerson2005,keller2010,bishop2010,quarteroni2017,piersanti,regazzoni2022,africa2023}, on space-dependent conduction velocities~\cite{usyk2002,goktepe2010,trayanova2011,sermesant2012}, or on the introduction of a fast endocardial layer~\cite{kotikanyadanam2010,baillargeon2014,neic2017,lee2019,gillette2021,gillette2022,verzicco2022,fedele2023,piersanti20223,bucelli2023,zappon2024}.
Other studies include an explicit representation of the His-Purkinje system, either rule-based~\cite{vergara2016,abboud1991,berenfeld1998,vigmond2007construction,ijiri2008,ten2008,romero2010,sebastian2012,bordas2012bidomain,behradfar2014,costabal2016,landajuela2018,strocchi2020,peirlinck2021,strocchi2022,reza2024assessing} or patient-specific~\cite{palamara_anti,camara2010,palamara_ortho,camps2024}. However, only a limited number of works have developed models that account for antidromic propagation, thus allowing a bidirectionally-coupled Purkinje-myocardium model. 

The main issue in modeling this interplay is capturing all the possible reentries of the electrical signal between the network and the muscle. Reentrant processes have been modeled using the Monodomain equation both for network and for muscular propagation~\cite{vergara2016,berenfeld1998,ten2008,vigmond2016,vigmond2021}, or using Bidomain models~\cite{behradfar2014,landajuela2018}. The work in~\cite{palamara_anti} proposed an algorithm for Purkinje-muscle coupling using an \textit{Eikonal-Eikonal} model~\cite{franzone3}, which accounts for antidromic propagation arising from possible muscular sources. However, this algorithm makes the simplifying assumption that the signal propagating from the network to the muscle cannot reenter the network, as would occur, for instance, in the presence of bundle branch blocks~\cite{berenfeld1998,ten2008,vigmond2016,vigmond2021}. Recently, a reaction-Eikonal model~\cite{neic2017}, that incorporates the Purkinje network through a monolithic approach~\cite{strocchi2020,strocchi2022}, yet without a detailed discussion of antidromic reentries, has been used to replicate clinical outcomes, such as for \ac{CRT}, His bundle pacing, and \ac{ECGs}. In conclusion, to the best of our knowledge, no studies have been conducted using Eikonal models to address the modeling of bidirectional propagation with a discussion on reentries. Models of this kind would be useful, particularly thanks to the computational efficiency of solvers for Eikonal equations.

In this work, we present a new Eikonal-Eikonal partitioned scheme for Purkinje-muscle electrophysiology simulations, supporting both orthodromic and antidromic propagation, with a particular focus on reentries. We rely on the Eikonal model for the network and on the Eikonal-diffusion model for the muscle. The main motivations underlying our choices lie in: (i) the modularity offered by partitioned schemes, compared to the monolithic approach, allowing to solve the two problems independently using distinct solvers; and (ii)~the computational advantage offered by the Eikonal models, which significantly reduce computational costs compared to the Monodomain and Bidomain ones. 

We remark that Monodomain and Bidomain models have the advantage of inherently handling bidirectional propagation and PMJ delays~\cite{vergara2016}, while an Eikonal-Eikonal model necessitates specific algorithms to explicitly address these interactions. In this regard, our novel algorithm improves the versatility and generality of the one proposed in~\cite{palamara_anti}, enabling reentries of the electrical signal between the network and the muscle. Moreover, our algorithm requires the development of a new \textit{pseudo-time method} for the Eikonal-diffusion problem, as classic solvers would not be accurate in this coupling context.

We showcase the accuracy and the efficiency of the proposed computational framework through numerical experiments in realistic geometries and under pathological or therapeutic conditions, such as bundle branch blocks and \ac{CRT}.

This paper is organized as follows. Section~\ref{sec:math} outlines the mathematical and numerical models employed in this work. Section~\ref{sec:setting} details the general parameters and computational setup. Section~\ref{sec:results} presents the results of numerical simulations in various pathological scenarios. Finally, Section~\ref{sec:conclusions} provides a summary of the conclusions.
	
	\section{Mathematical and numerical models}
	\label{sec:math}
	In this section, we introduce the mathematical models used to describe the propagation of electrical signals within the cardiac conduction system. In the subsequent sections we examine in detail the following modeling aspects:
\begin{itemize}
	\item the model for the activation of the His-Purkinje system through Eikonal equation;
	\item the model for the activation of the myocardium using the Eikonal-diffusion equation;
	\item the Purkinje-myocardium coupling to manage reentries between the two domains, resulting from orthodromic and antidromic propagation, employing a novel partitioned scheme for the bidirectional coupling.
\end{itemize}

\subsection{Electrical activation of the Purkinje}
\subsubsection{Continuous problem}
We use the Eikonal equation~\cite{palamara_anti,vergara2016} to describe the propagation of electrical signals within the Purkinje fibers, represented as a one-dimensional domain $\Omega_\mathrm{p} \subset \mathbb{R}^3$ (see Figure~\ref{fig:domain}(right)):
\newline
\\
$\text{\quad \quad \quad} \text{Given } \Gamma_0 \subset \partial\Omega_\mathrm{p} \text{ and } u_\mathrm{p}^0:\Gamma_0 \rightarrow [0,T] ;\text{ find  } u_\mathrm{p}(\vb*{x}):\Omega_\mathrm{p} \rightarrow [0,T] \text{  such that}$
\begin{equation}
	\label{eq:eikonal1D}
	\begin{dcases}
		\begin{aligned}
			c_\mathrm{p} \bigg|\frac{\partial u_\mathrm{p}}{\partial s} \bigg| = 1 &\quad \quad \text{in } \Omega_\mathrm{p}, \\
			u_\mathrm{p} = u_\mathrm{p}^0 &\quad \quad \text{on } \Gamma_0, 
		\end{aligned}
	\end{dcases}
\end{equation}
where $u_\mathrm{p}$ represents the unknown activation times of the network, $c_\mathrm{p}$ is the conduction velocity, and $s$ is the local curvilinear coordinate along the network. $\Gamma_0$ is the set of points generating the front. In healthy scenarios, $\Gamma_0$ coincides with the atrioventricular node (AV node), whereas in general it may be composed by PMJs supporting antidromic propagation.

Equation~\eqref{eq:eikonal1D} provides the activation times at any point of the Purkinje network, including the \ac{PMJs}. Conduction blocks are modeled as disconnections in the network.

\subsubsection{Numerical discretization: Fast Marching method}
\label{subsubsec:discretization_FMM}
We solve the 1D Eikonal problem using the \ac{FMM}, proposed for the first time in~\cite{sethian} to solve the pure Eikonal equation and already applied to the context of cardiac electrophysiology~\cite{vergara2016,pernod2011,grandits2020}. In the 1D setting, the \ac{FMM} reduces to the Dijkstra's shortest path algorithm~\cite{dijkstra}.

Notably, this algorithm automatically updates the activation time for those nodes in the network where a condition is prescribed (i.e. nodes belonging to $\Gamma_0$), but a signal coming from another stimulation point arrives earlier. In other words, if a stimulus is applied to a point where the solution to~\eqref{eq:eikonal1D} would result in an activation earlier than the stimulus itself, then that stimulus is disregarded by the algorithm. 

\subsection{Electrical activation of the myocardium}
\subsubsection{Continuous problem}
\begin{figure}
	\centering
	\includegraphics[width=1\textwidth]{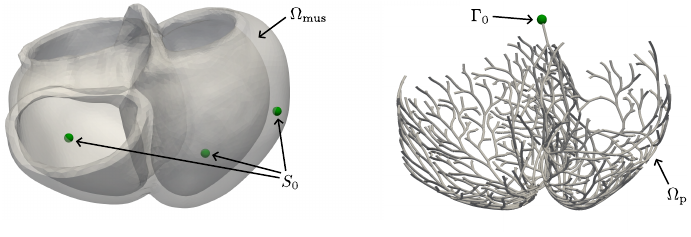}
	\caption{Ventricular myocardial domain $\Omega_\mathrm{mus}$ (left) with three sources at $S_0$ and corresponding Purkinje domain $\Omega_\mathrm{p}$ (right) with a source at $\Gamma_0$.}
	\label{fig:domain}
\end{figure}

The propagation of the electrical signal in the myocardium $\Omega_\mathrm{mus} \subset \mathbb{R}^3$, see Figure~\ref{fig:domain}(left), is modeled by the Eikonal-diffusion equation~\cite{quarteroni2017,franzone3}. This model is used to reconstruct the macroscopic propagation of action potential excitation wavefronts during the depolarization phase, allowing to compute the activation time as a spatial function within the myocardium. 

We define the activation time $u_\mathrm{m}=u_\mathrm{m}(\vb*{x})$ as the instant at which the depolarization wavefront reaches the point $\vb*{x} \in \Omega_\mathrm{mus}$. The Eikonal-diffusion model reads:
\newline
\\
$\text{\quad \quad \quad} \text{Given } u_\mathrm{m}^0 \in [0,T] \text{ for some } \vb*{x} \in S_0; \text{ find  } u_\mathrm{m}(\vb*{x}):\Omega_\mathrm{mus} \rightarrow [0,T] \text{  such that}$
\begin{equation}
	\label{eq:eikonal3D}
	\begin{dcases}
		\begin{aligned}
			c_\mathrm{f} \sqrt{\nabla u_\mathrm{m} \cdot \vb{\Sigma} \nabla u_\mathrm{m}} - \nabla \cdot (\vb{\Sigma} \nabla u_\mathrm{m}) = 1 &\quad \quad \text{in } \Omega_\mathrm{mus}, 
			\\
			(\vb{\Sigma} \nabla u_\mathrm{m}) \cdot \vb{n} = 0 &\quad \quad \text{on } \partial\Omega_\mathrm{mus} \text{\textbackslash} S_0, 
			\\
			u_\mathrm{m} = u_\mathrm{m}^0 &\quad \quad \text{on } S_0.
		\end{aligned}
	\end{dcases}
\end{equation}
Here, $\vb{n}$ is the outward directed unit vector normal to the boundary $\partial\Omega_\mathrm{mus}$ of the domain $\Omega_\mathrm{mus}$, while $S_0$~is a portion of the physical boundary where the activation time is prescribed. Under normal propagation, $S_0$~coincides with the \ac{PMJs} responsible for orthodromic propagation, because the signal traveling along the Purkinje network enters the muscle through these junctions. In cases of dysfunctions, additional muscular sources could be present, such as ectopic stimulation sites~\cite{palamara_anti}. In~\eqref{eq:eikonal3D}, $\vb{\Sigma} \in \mathbb{R}^{3 \times 3}$ is the anisotropic tensor accounting for the orientation of muscular myofibers, and it is related to the conductivity tensor $\vb{D}(\vb*{x})$ through the membrane capacitance $C_\mathrm{m}$ and the surface-to-volume ratio $\chi_\mathrm{m}$:
\begin{equation}
	\label{eq:conductivity}
	\vb{\Sigma} = \frac{\vb{D}}{\chi_\mathrm{m} C_\mathrm{m}} \text{\quad with \quad} \vb{D}(\vb*{x}) = \sigma_\mathrm{f} \vb*{f}_0 \otimes \vb*{f}_0 + \sigma_\mathrm{s} \vb*{s}_0 \otimes \vb*{s}_0 + \sigma_\mathrm{n} \vb*{n}_0 \otimes \vb*{n}_0,
\end{equation}
where $\sigma_\mathrm{f},\sigma_\mathrm{s},\sigma_\mathrm{n}$ are the conductivities along the fibers (longitudinal) $\vb*{f}_0$, the sheets (transversal) $\vb*{s}_0$, and the cross-fibers (normal) $\vb*{n}_0$ directions, respectively. Accurately modeling myocardial fiber arrangement is essential to reproduce the anisotropic propagation occurring in the muscle~\cite{piersanti,bayer2012novel}. We prescribed muscular fibers, employing the strategy proposed by~\cite{piersanti,lifex-fiber}, through the algorithm by Bayer et al.~\cite{bayer2012novel}, which is a Laplace-Dirichlet rule-based method. Finally, the parameter $c_\mathrm{f}$,
uniform over the domain $\Omega_\mathrm{mus}$, is the velocity of the action potential depolarization planar wavefront along the fiber direction in an infinite cable, under the assumption of a unit surface-to-volume ratio, membrane capacitance, and conductivity~\cite{franzone1,franzone2} and normalized with respect to the diffusion parameter.

The first term in the first equation of~\eqref{eq:eikonal3D} is the classical Eikonal term, describing the propagation of a front in an anisotropic domain (where the anisotropy is captured by $\vb{D}$). The second term is a generalized Laplacian representing anisotropic diffusion, establishing a connection between the speed of the wavefront and its curvature. Propagation is faster when the wavefront is concave, and slower when it is convex~\cite{quarteroni103,franzone2}.

Problem~\eqref{eq:eikonal3D} is very convenient from a computational standpoint, if compared to the Monodomain and the Bidomain models. Firstly, it requires solving a single steady PDE, albeit non-linear, as opposed to the Mondomain and Bidomain problems, which are time-dependent and must be coupled to ODE systems for ionic modeling~\cite{quarteroni103,franzone3}. More importantly, the activation time, unlike the transmembrane potential, typically does not feature sharp fronts or strong gradients. Therefore, the Eikonal model does not demand for a fine spatial resolution~\cite{colli2014,strocchi,stella2022,gander2023}, enabling the simulation of the activation of large volumes of cardiac tissue with relatively low computational costs.

\subsubsection{Numerical discretization: a novel pseudo-time method}
\label{subsec:novel_psuedo}
The numerical solution of~\eqref{eq:eikonal3D} requires a discretization with a proper treatment of the non-linearity. Classical non-linear solvers, such as Newton's method, typically have convergence issues in this context, due to the difficulty of finding an appropriate initial guess to problem~\eqref{eq:eikonal3D}~\cite{franzone3}. To circumvent this issue, one can numerically solve~\eqref{eq:eikonal3D} by recovering the steady-state solution of the following parabolic pseudo-time problem~\cite{franzone3,quarteroni2006}:
\begin{equation}
	\label{eq:pseudo_time}
	\frac{\partial u_\mathrm{m}}{\partial t} + c_\mathrm{f} \sqrt{\nabla u_\mathrm{m} \cdot \vb{\Sigma} \nabla u_\mathrm{m}} - \nabla \cdot (\vb{\Sigma} \nabla u_\mathrm{m}) = 1 \quad \quad \text{in } \Omega_\mathrm{mus}. 
\end{equation}
This is known as \textit{pseudo-time} method. In~\eqref{eq:eikonal3D_pseudo_time}, we propose a different version of the pseudo-time method, introduced to properly handle the Purkinje-muscle coupled propagation, since standard versions would not be accurate in the context of our new coupling algorithm, as detailed in what follows.

Starting from Equation~\eqref{eq:pseudo_time}, the time discretization employing a fully implicit \ac{BDF} of order $\sigma$ (say BDF$\sigma$; $\sigma=1,2,...$), along with proper boundary and initial conditions, reads:
\newline
\\
$\text{\quad \quad \quad} \text{At each psuedo-time } t^n \text{, given } S^{n+1} \subset S_0 \text{ and } u_\mathrm{m}^0: S_0 \rightarrow \mathbb{R}; \text{ find  } u^{n+1} \text{  such that}$
\begin{equation}
	\label{eq:eikonal3D_pseudo_time}
	\begin{dcases}
		\begin{aligned}
			\frac{\alpha_\mathrm{BDF} u^{n+1} - u^{n}_\mathrm{BDF\sigma}}{\Delta t} + c_\mathrm{f} \sqrt{\nabla u^{n+1} \cdot \vb{\Sigma} \nabla u^{n+1}} - \nabla \cdot (\vb{\Sigma} \nabla u^{n+1}) & = 1 \quad \quad \hspace{0.1cm}\text{in } \Omega_\mathrm{mus}, 
			\\
			(\vb{\Sigma} \nabla u^{n+1}) \cdot \vb{n} & = 0 \quad \quad \hspace{0.1cm}\text{on } \partial\Omega_\mathrm{mus} \text{\textbackslash} S^{n+1}, 
			\\
			u^{n+1} & = u_\mathrm{m}^0 \hspace{-0.2cm} \quad \quad \hspace{0.1cm}\text{on } S^{n+1},
		\end{aligned}
	\end{dcases}
\end{equation}
where $u^n\simeq u_\mathrm{m}(\vb*{x},t^n)$ is the numerical solution, $\Delta t = t^{n+1} - t^n$, while $\alpha_\mathrm{BDF}$ and $u^{n}_\mathrm{BDF\sigma}$ are the terms associated to the BDF$\sigma$ time discretization\footnote{$\alpha_\mathrm{BDF}$ is a constant coefficient, while $u^{n}_\mathrm{BDF\sigma}$ refers to the value of the function $u$ at the current, or possibly earlier, time steps, or even a combination of them. For instance, for BDF1, $\alpha_\mathrm{BDF} = 1$ and $u^{n}_\mathrm{BDF\sigma} = u^n$, while for BDF2, $\alpha_\mathrm{BDF} = 3$ and $u^{n}_\mathrm{BDF\sigma} = 4u^n - u^{n-1}$, and so on.}~\cite{quarteroni2010numerical,forti2015semi}. The novelty of our approach lies in the location where Dirichlet-type conditions are imposed. Specifically, the collection of all $N$ stimulation points $S_0 = \{\textbf{x}_i\}_{i=1}^N$, is replaced now by $S^{n+1}$, containing only the subset of $S_0$ related to \textit{active} stimuli. At time $t^{n+1}$, active stimuli are identified by:
\begin{equation}
	\mathbf{x}_i \in S^{n+1} \iff u_\text{m}^0(\mathbf{x}_i) < t^{n+1} \text{\quad and \quad} u_\text{m}^0(\mathbf{x}_i) < u^n(\mathbf{x}_i)\;.
	\label{eq:active_condition}
\end{equation}
Dirichlet-type conditions are imposed only on the stimulation points that the algorithm identifies as active, whereas inactive stimuli are disregarded. Therefore, we underline that the solution of \eqref{eq:eikonal3D_pseudo_time} is deliberately not consistent with the continuous problem \eqref{eq:eikonal3D}, since the steady-state solution only satisfies the condition $u_\text{m} = u_\text{m}^0$ on a subset of $S_0$. Thus, \eqref{eq:eikonal3D_pseudo_time} can be seen both as an alternative approach to Eikonal modeling of cardiac activation and as a solving methodology. A rigorous derivation of the pseudo-time-continuous counterpart to~\eqref{eq:eikonal3D_pseudo_time} is not addressed in this paper.

\begin{figure}
	\centering
	\includegraphics[width=1\linewidth]{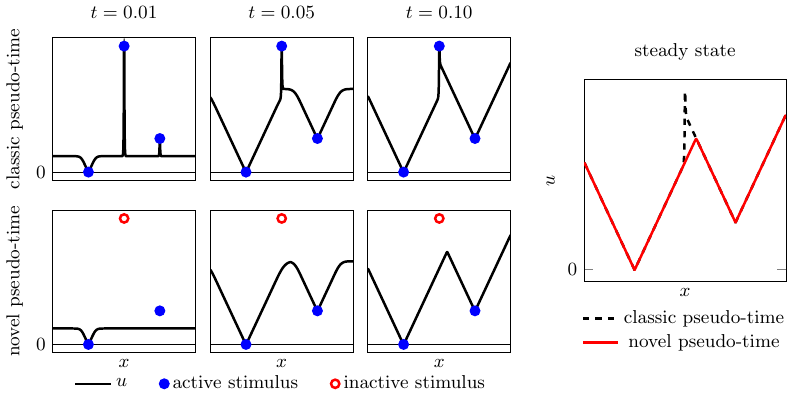}
	\caption{Numerical solution returned by the novel pseudo-time method in a one-dimensional example, compared to the classic pseudo-time method. The spatial domain is $x \in [0,L]$. Three stimuli are imposed: two are classified as active stimuli by the novel algorithm (in blue) and one as inactive (in red). At each pseudo-time step the solution is updated and stimulation points are identified as either active or inactive according to~\eqref{eq:active_condition}. Notably, the steady-state solution obtained with the classic pseudo-time method is physically meaningless, because the latest stimulus should be disregarded, as its activation time is later than the one induced by the other two.}
	\label{fig:new_pseudo_time}
\end{figure}

Figure~\ref{fig:new_pseudo_time} helps clarifying the rationale behind the development of our novel solver. The figure compares the evolution of the solutions obtained by the classic and novel pseudo-time methods for the same set of stimuli. The novel algorithm correctly identifies one stimulus as being inactive, whereas the classic algorithm would lead to a non-physical activation delay in correspondence of that stimulus. The novel pseudo-time formulation is essential for obtaining physically meaningful solutions during the coupling with the Purkinje network. 

For the linearization and space discretization of~\eqref{eq:eikonal3D_pseudo_time} we use the standard Newton method and the Galerkin-Finite Element method~\cite{piersanti2024defining}.

\subsection{Purkinje-muscle fully coupled problem}
\label{subsec:coupling}
\begin{figure}
	\centering
	\includegraphics[width=1\textwidth]{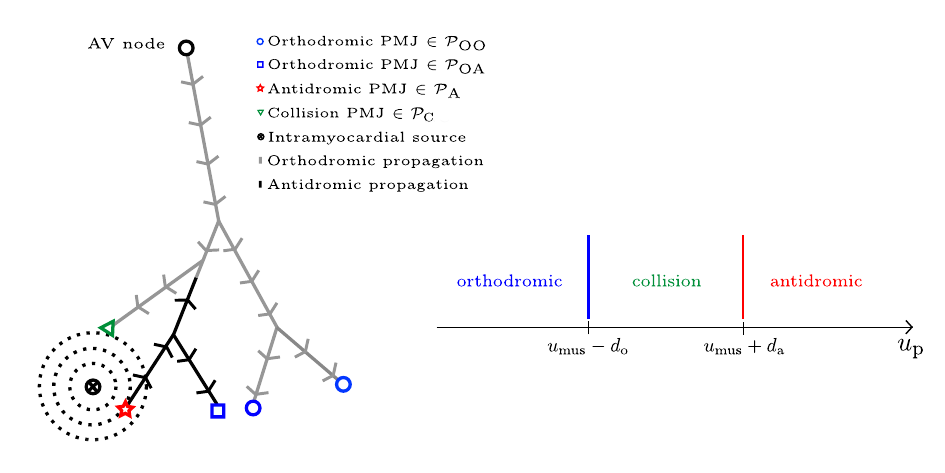}
	\caption{\acf{PMJs} classification in presence of opposite wavefronts propagating between the Purkinje network and the myocardium. Right: graphical representation of the function {\textbf{classify\_pmjs}}(\textit{pmjs}\texttt{,}$u_{\mathrm{p}}$\texttt{,}$u_{\mathrm{m}}$), defined in Algorithm~\ref{alg:pmj_classification}. In the notation, AV node stands for atrioventricular node. Two types of orthodromic PMJs are introduced: $\mathcal{P}_{\text{OO}}$ is the collection of orthodromic \ac{PMJs} activated by an electrical signal coming from the AV node, while $\mathcal{P}_{\text{OA}}$ is the collection of orthodromic \ac{PMJs} activated by a signal coming from an antidromic PMJ.}
	\label{fig:antidromic}
\end{figure}

To address the coupling between the problems~\eqref{eq:eikonal1D} and~\eqref{eq:eikonal3D}, we introduce the following functionals, which map the stimuli to the corresponding activation times both in the myocardium and in the Purkinje network:
\begin{equation}
	\nonumber
	\small
	\begin{aligned}
		\dutchcal{E}: S_0 \longmapsto u_\mathrm{m} \quad \quad \text{\quad \quad} \dutchcal{P}: \Gamma_0 \longmapsto u_\mathrm{p} .
	\end{aligned}
\end{equation}
Our starting point is the coupling strategy presented in~\cite{palamara_anti}. The interface conditions of the coupled problem are related to the continuity of the activation times $u_\mathrm{p}$ and $u_\mathrm{m}$ at the PMJs. The latter are classified into:

\begin{itemize}
	\item \textit{antidromic \ac{PMJs}}, which are activated by the front coming from the muscle, and delay the signal by~$d_\mathrm{a}$ before transmitting it to the Purkinje network.
	\item \textit{orthodromic \ac{PMJs}}, which are activated by the front coming from the AV node or emerging in the network due to antidromic propagation, and delay the signal by $d_\mathrm{o}$ before transmitting it to the myocardium.
	\item \textit{collision \ac{PMJs}}\footnote{This is a new category, introduced in this work, to speed up and correctly classify PMJs in Algorithm~\ref{alg:pseudocode_bis}.}, which arise when the muscular and the network wavefronts collide at the PMJ itself, so they are neither orthodromic nor antidromic. This is possible due to the presence of a delay at the \ac{PMJs}.
\end{itemize} 
Figure~\ref{fig:antidromic} schematically illustrates the different types of \ac{PMJs} resulting from orthodromic and antidromic propagation. 

We propose hereafter Algorithm~\ref{alg:pseudocode_bis}, designed to overcome the limitations of \cite{palamara_anti}. The latter made the assumption that the antidromic activation of the \ac{PMJs} is induced exclusively by the muscular sources. This implies that the muscular pathways that activate the network could only be generated by the muscular sources. In other words, with the algorithm of~\cite{palamara_anti}, a signal originating from the AV node and entering the muscle through a PMJ cannot reenter the network. Conversely, our approach is able to account for potential occurrence of reentries.  Algorithm~\ref{alg:pseudocode_bis} allows for an excitation front to originate at the AV node, enter the muscle through a PMJ and subsequently reenter the Purkinje network at a different PMJ. This is achieved by incorporating subiterations in the algorithm, thus repeatedly solving the two Eikonal problems and exchanging information at the \ac{PMJs} until the maximum number of iterations is reached, and it is made possible by the use of the new pseudo-time method described in Section~\ref{subsec:novel_psuedo}. We fix the number of iterations to $N=N_{\max}$, although this can be straightforwardly improved through a suitable stopping criterion.

The \ac{PMJs} classification, showed in Figure~\ref{fig:antidromic}, is updated through the function \textbf{classify\_pmjs}, illustrated in Algorithm~\ref{alg:pmj_classification}.
\begin{algorithm}
	\caption{1D Eikonal - 3D Eikonal strong coupling}
	\label{alg:pseudocode_bis}
	\begin{algorithmic}[1]
		\small
		\vspace{0.1cm}
		\State{Notation: $AVN$ = atrioventricular node, $\mathcal{M}$ = collection of the possible muscular sources, $pmj_{\mathrm{o}_i}$ =  $i$-th orthodromic PMJ, $pmj_{\mathrm{a}_i}$ =  $i$-th antidromic PMJ, $\texttt{u}_{\texttt{m}}$ = activation time in the muscle, $\texttt{u}_{\texttt{p}}$ = activation time in the His-Purkinje system, $\texttt{pmj}_{\texttt{o}}$/$\texttt{pmj}_{\texttt{a}}$/$\texttt{pmj}_{\texttt{c}}$ = orthodromic/antidromic/collsion PMJs}
		\vspace{0.3cm}
		\State{\textbf{Input:} $AVN, \mathcal{M}, N_{\max}$}
		\vspace{0.1cm}
		\State{\textbf{Output:} $\texttt{u}_{\texttt{p}}, \texttt{u}_{\texttt{m}}, \texttt{pmjs}$}
		\vspace{0.3cm}
		\State{$\texttt{u}_{\texttt{p}} \gets \mathdutchcal{P}~(\{AVN\})$}
		\vspace{0.2cm}
		\State{$\texttt{pmjs} \gets$ get Purkinje-muscle junctions (PMJs)}
		\vspace{0.2cm}
		\State{{\textbf{for}} (\texttt{$i<N_{\max}$}) {\textbf{do}}}
		\vspace{0.2cm}
		\State{$\quad \texttt{u}_{\texttt{m}} \gets \mathdutchcal{E}~(\mathcal{M}\cup\{pmj_{\mathrm{o}_1}, pmj_{\mathrm{o}_2}, ... , pmj_{\mathrm{o}_N} \})$}
		\vspace{0.2cm}
		
		\State{$\quad \texttt{[pmj}_{\texttt{o}} \texttt{,pmj}_{\texttt{a}} \texttt{,pmj}_{\texttt{c}} \texttt{]} \gets {{\textbf{classify\_pmjs}}} \texttt{(pmjs,u}_{\texttt{p}} \texttt{,u}_{\texttt{m}} \texttt{)}$}
		\vspace{0.2cm}
		
		\State{$\quad \texttt{u}_{\texttt{p}} \gets \mathdutchcal{P}~(\{AVN, pmj_{\mathrm{a}_1}, pmj_{\mathrm{a}_2}, ... , pmj_{\mathrm{a}_M} \})$}
		\vspace{0.2cm}
		
		\State{$\quad \texttt{[pmj}_{\texttt{o}} \texttt{,pmj}_{\texttt{a}} \texttt{,pmj}_{\texttt{c}} \texttt{]} \gets {{\textbf{classify\_pmjs}}} \texttt{(pmjs,u}_{\texttt{p}} \texttt{,u}_{\texttt{m}} \texttt{)}$}
		
		\vspace{0.2cm}
		\State{{\textbf{end for}}}
	\end{algorithmic}
\end{algorithm} 
\begin{algorithm}
	\caption{Function ${{\textbf{classify\_pmjs}}} \texttt{(pmjs,u}_{\texttt{p}} \texttt{,u}_{\texttt{m}} \texttt{)}$}
	\label{alg:pmj_classification}
	\begin{algorithmic}[1]
		\small
		\vspace{0.1cm}
		\State{{\textbf{for}} (\texttt{pmj:pmjs}) {\textbf{do}}}
		\vspace{0.2cm}
		\State{{\textbf{\quad if}} $(\texttt{u}_{\texttt{p}} \texttt{[pmj] >= u}_{\texttt{m}} \texttt{[pmj] + antidromic\_delay)}$ {\textbf{then}}}
		\vspace{0.1cm}
		\State{\quad \quad $\texttt{PmjType[pmj] = "Antidromic"}$}
		\vspace{0.2cm}
		
		\State{{\textbf{\quad else if}} $(\texttt{u}_{\texttt{p}} \texttt{[pmj] <= u}_{\texttt{m}} \texttt{[pmj] - orthodromic\_delay})$ {\textbf{then}}}
		\vspace{0.1cm}
		\State{\quad \quad $\texttt{PmjType[pmj] = "Orthodromic"}$}
		\vspace{0.2cm}
		
		\State{{\textbf{\quad else}}}
		\vspace{0.1cm}
		\State{\quad \quad $\texttt{PmjType[pmj] = "Collision"}$}
		\vspace{0.1cm}
		\State{{\textbf{\quad end if}}}
		\vspace{0.1cm}
		\State{{\textbf{end for}}}
	\end{algorithmic}	
\end{algorithm}
\clearpage
	
	\section{Simulations setting}
	\label{sec:setting}
	In this section, we briefly introduce the general setting of the numerical simulations that will be presented in Section~\ref{sec:results}. Section~\ref{subsec:domains} describes the computational domains considered for the 3D ventricular geometry and the 1D Purkinje network, while Section~\ref{subsec:params} summarizes the common parameter settings for all the numerical simulations. Finally, Section~\ref{subsec:code} is dedicated to the computational environment adopted to implement the algorithms and to perform the simulations.

\subsection{Computational domains}
\label{subsec:domains}
\begin{figure}
	\centering
	\includegraphics[width=1\textwidth]{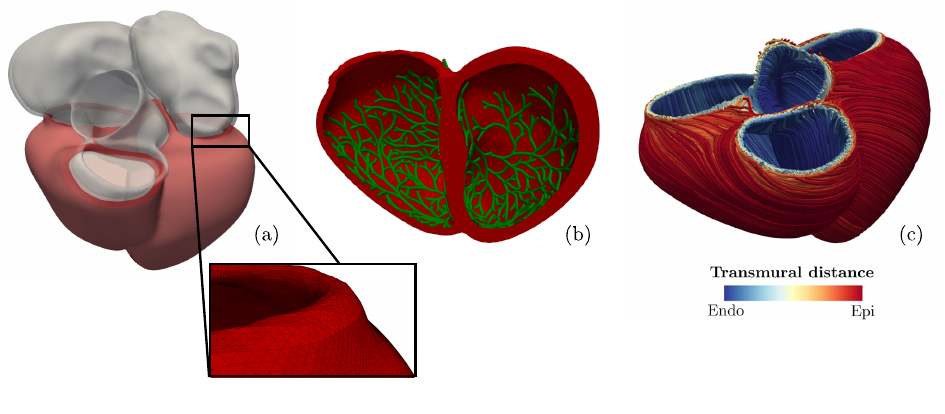}
	\caption{(a) Real biventricular geometry reconstructed from four-chamber heart of a patient~\cite{strocchi}, with tetrahedral mesh discretization. (b) A section of the ventricles with the generated Purkinje network. (c) Muscular fibers visualization. The fiber field is visualized as streamlines, colored according to the transmural distance, going from the endocardium (Endo, in blue) to the epicardium (Epi, in red).}
	\label{fig:geometry}
\end{figure}
The 3D cardiac computational domain, depicted in Figure~\ref{fig:geometry}, is a patient-specific human biventricular geometry, reconstructed from a \ac{CT} scan of a patient. This geometry is obtained from a publicly available cohort of four-chamber hearts presented in~\cite{strocchi}. The mesh was segmented from a \ac{CT} scan of the heart of an 83-years-old male who was recruited for a \ac{CRT} upgrade following heart failure. Since this is not a healthy heart, it exhibits certain non-physiological features such as muscular hypertrophy and myocardial thickening. Moreover, endocardial structures, such as papillary muscles, are not included in the segmentation.

The extraction of the ventricular portion and the volumetric mesh generation (with an average edge length of \SI{0.8}{\milli\meter}) was performed using \texttt{vmtk}\footnote{\protect\url{http://www.vmtk.org/}} and the methods described in~\cite{fedele2021polygonal}. The main characteristics of the obtained mesh are summarized in Table~\ref{table:mesh}.

Incorporating a patient-specific His-Purkinje system represents a distinct challenge, since currently no imaging technique is able to reconstruct a patient-specific network for human hearts. We used the open-source \texttt{Python} library \texttt{fractal-tree}\footnote{\protect\url{https://github.com/fsahli/fractal-tree}} to generate the 1D Purkinje network~\cite{costabal2016}. This approach employs a semi-automatic algorithm based on a fractal law, designed to generate a network standing on an endocardial surface discretized by triangles, mimicking the tree structure of the Purkinje itself. An example of generated network is reported in Figure~\ref{fig:geometry}(b). The absence of endocardial structures in the considered 3D geometries prevents an accurate reconstruction of the network, allowing only to capture its principal features. Some significant elements, particularly the moderator band in the right ventricle, have been omitted.

\subsection{Simulation parameters}
\label{subsec:params}
\begin{table}
	\centering
	\begin{tabular}{ c  c }
		\hline
		\textbf{element shape} & tetrahedra \\
		\textbf{\#elements} & 5'268'810 \\
		\textbf{\#nodes(DoFs)} & 927'097 \\
		\textbf{aspect ratio} & 1.57 $\pm$ 0.53 \\ 
		\textbf{scaled jacobian} & 0.63 $\pm$ 0.14 \\ 
		\textbf{average cell diameter} & \SI{0.927}{\milli\meter}\\
		\textbf{min cell diameter} & \SI{0.377}{\milli\meter} \\
		\textbf{max cell diameter} & \SI{1.365}{\milli\meter} \\ 
		\hline
	\end{tabular}
	\caption{Characteristics of the 3D mesh used in numerical simulations. Aspect ratio and scaled jacobian are two mesh quality metrics, defined as in~\cite{mesh_quality}, indicating that the mesh has acceptable regularity.}
	\label{table:mesh}
\end{table}
\begin{table}
	\centering
	\begin{tabular}{c c c}
		\hline
		\textbf{Description} & \textbf{Parameter} & \textbf{Reference} \\
		\hline
		\\[-0.3cm]
		Stimulation sites and shape & & \\
		Conduction block sites & & \\
		& $\sigma_\mathrm{f} =$ \SI[per-mode=symbol]{1.00e-4}{\square\meter\per\second} & \\
		Conductivity values &
		$\sigma_\mathrm{s} =$ \SI[per-mode=symbol]{0.44e-4}{\square\meter\per\second} &      \cite{piersanti} \\
		& $\sigma_\mathrm{n} =$ \SI[per-mode=symbol]{0.11e-4}{\square\meter\per\second} & \\
		Depolarization velocity of the myocardium & $c_\mathrm{f} = \SI{60}{\per\second^{-1/2}}$ & \cite{piersanti} \\
		Conduction velocity of the Purkinje & $c_\mathrm{p} = \SI[per-mode=symbol]{4}{\meter\per\second}$ & \cite{cavero2023,vergara2016} \\
		Orthodromic PMJ delay & $d_\mathrm{o} = \SI{10e-3}{\second}$ & \cite{vergara2016} \\
		Antidromic PMJ delay & $d_\mathrm{a} = \SI{2e-3}{\second}$ & \cite{vergara2016} \\
		\hline
	\end{tabular}
	\caption{Summary of the main physical parameter used in the numerical simulations.}
	\label{table:parameters}
\end{table}
Table~\ref{table:parameters} reports the main physical parameters useful for solving the Eikonal model. Following~\cite{piersanti}, we select the conductivity values in the myocardium to obtain the conduction velocities of \SI[per-mode=symbol]{0.6}{\meter\per\second} in the fiber direction, \SI[per-mode=symbol]{0.4}{\meter\per\second} in the sheet direction and \SI[per-mode=symbol]{0.2}{\meter\per\second} in the normal direction~\cite{piersanti,cavero2023}. We selected the conduction velocity in the Purkinje fibers to be \SI[per-mode=symbol]{4}{\meter\per\second}~\cite{verzicco2022,cavero2023,palamara_ortho}. At the \ac{PMJs}, the signal experiences a delay of $\SIrange[range-phrase=-]{10}{15}{\milli\second}$ for orthodromic propagation and a shorter delay of $\SIrange[range-phrase=-]{2}{3}{\milli\second}$ for antidromic propagation~\cite{vergara2016}.

Finally, pointwise stimuli in the muscle are applied at the \ac{PMJs}, while spherical stimuli of radius \SI{0.5}{\milli\meter} are imposed in correspondence of possible ectopic or artifical stimuli. The placement of muscular sources and possible conduction blocks is evaluated case-by-case and passed as input through their point coordinates.

\subsection{Simulation environment}
\label{subsec:code}
The numerical framework presented in Section~\ref{sec:math} has been implemented in $\lifex$~\cite{africa,lifex-ep,bucelli2024lifex}\footnote{\url{https://lifex.gitlab.io}}, a high-performance \texttt{C++} finite element library for cardiovascular modeling, based on \texttt{deal.II}~\cite{dealii-original,dealii-requested,dealii-version9.4}.

Simulations have been performed using the HPC facilities of the Mathematics Department at Politecnico di Milano, running on one node with 20 cores 4x Xeon E5-2640 v4 (2.4GHz) with 384GB of RAM. Running $N_{\max}=3$ iterations\footnote{The choice $N_{\max}=3$, in Algorithm~\ref{alg:pseudocode_bis}, is done to balance reliable results with reduced computational time.} of Algorithm~\ref{alg:pmj_classification} takes less than 2 $\mathrm{hours}$.
	
	\section{Numerical results}
	\label{sec:results}
\begin{figure}
	\centering
	\includegraphics[width=0.96\textwidth]{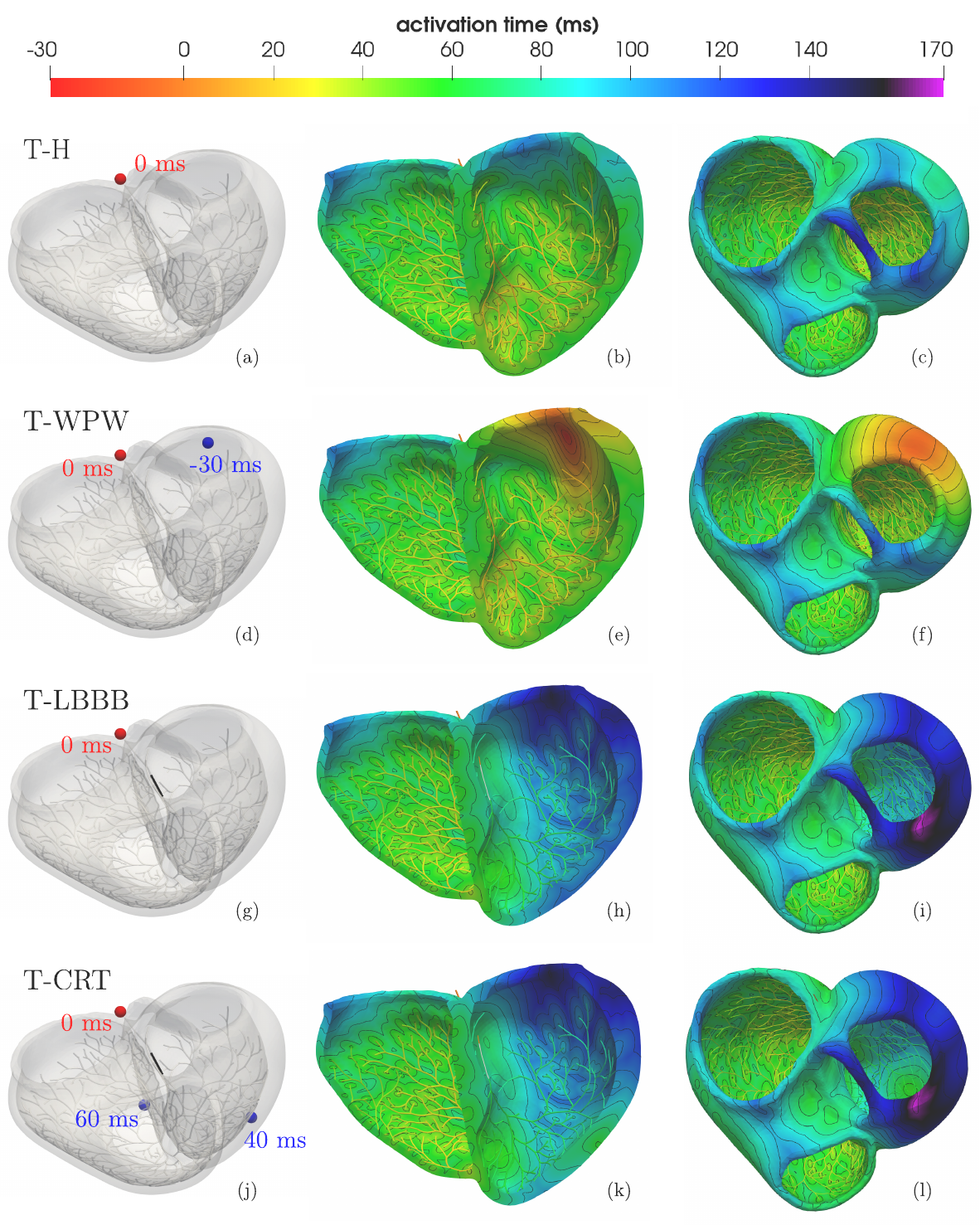}
	\caption{Comparative view of the four scenarios: T-H, T-WPW, T-LBBB, T-CRT. First column: simulation setup; second column: activation map of a section of the biventricular geometry; third column: activation map from top view. Isochrones are spaced~\SI{10}{\milli\second} apart. The bundle branch block is highlighted in black in figures (g) and (j). The activation time of the  atrioventricular node (AV node) is taken as reference starting at \SI{0}{\milli\second}.}
	\label{fig:result1}
\end{figure}

\begin{figure}
	\centering
	\includegraphics[width=1\textwidth]{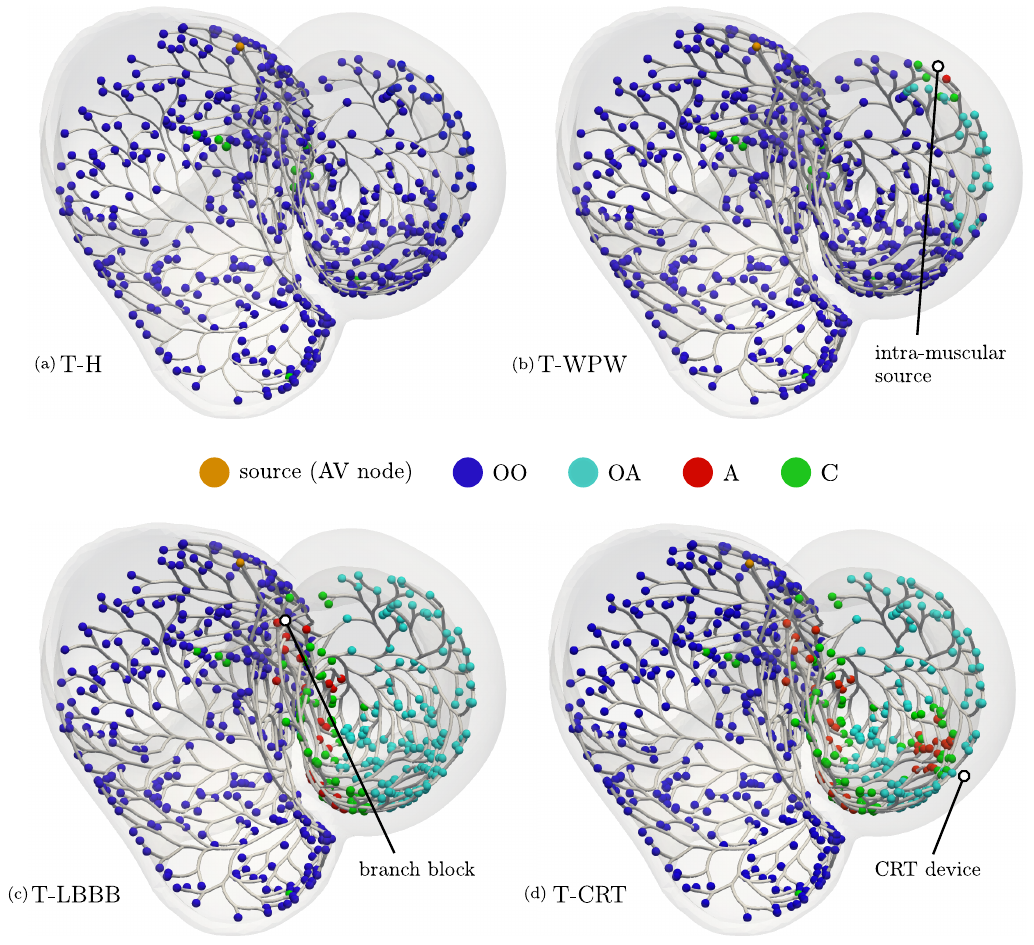}
	\caption{Graphical visualization of the \acf{PMJs} classification in the four test cases. \ac{PMJs} are represented as spheres of different colors according to their type: source (orange), orthodromic activated by AV node (OO, blue), orthodromic activated by antidromic propagation (OA, cyan), antidromic (A, red), collision (C, green). Non-orthodromic \ac{PMJs} are localized close to anomalous sources.}
	\label{fig:result4}
\end{figure}

	This section is dedicated to present numerical results from electrophysiology simulations using the model framework introduced in Section~\ref{sec:math} and based on the setting presented in Section~\ref{sec:setting}. Computational simulations offer a valuable tool for medical decision-making. This ability to simulate different scenarios with relative computational ease, without intervening on the patient, represents, in perspective, a notable application~\cite{quarteroni103,quarteroni2019mathematical,neic2017,strocchi2022,lee2018}. We consider four test (T) scenarios, illustrating both healthy and pathological propagation:
\begin{itemize}
	\item T-H: healthy propagation. The atrioventricular node represents the unique stimulation point for ventricular activation and only orthodromic propagation is present~\cite{fisiologia_ita}.
	\item T-WPW: Wolff-Parkinson-White syndrome. A pathological scenario characterized by an intramuscular source in addition to His-Purkinje activation~\cite{romano2010,palamara_anti,vergara2016}.
	\item T-LBBB: complete left bundle branch block. A pathological scenario characterized by a complete conduction block in the left bundle branch~\cite{romano2010,strocchi2022}.
	\item T-CRT: cardiac resynchronization therapy. A model therapeutical scenario, which combines both features of the previous test cases, introducing multiple intramuscular sources and a complete left bundle branch block~\cite{leyva2014,strocchi2022}.
\end{itemize}
Figures~\ref{fig:result1} and~\ref{fig:result4} compare the activation times and the \ac{PMJs} classification of the four test cases, respectively. The two figures are completed by Tables~\ref{table:postpro2} and~\ref{table:result5}, which collect synthetical indices related to the activation maps and to the \ac{PMJs} classification.

\begin{table}[]
	\centering
	\begin{tabular}{ c | c c c}
		\hline
		& $\mu \pm \sigma~(\mathrm{ms})$ & TAT & EAT \\ \hline
		\\[-0.3cm]
		T-H & $64.0 \pm 14.8$ & 140.7 & 27.6 \\
		T-WPW & $56.9 \pm 17.3$ & 120.5 & -30.0 \\
		T-LBBB & $86.0 \pm 26.4$ & 174.0 & 31.1 \\
		T-CRT & $83.8 \pm 25.6$ & 173.4 & 31.1 \\ 
		\hline
	\end{tabular}
	\caption{Average activation time ($\mu$), standard deviation ($\sigma$), total activation time (TAT) and earliest activation time (EAT) of the myocardium in the four test scenarios. All the times refer to $t=$~\SI{0}{\milli\second}, when the signal starts from the atrioventricular node.}
	\label{table:postpro2}
\end{table}
\begin{table}[]
	\centering
	\begin{tabular}{c | c c c c}
		\hline
		& T-H & T-WPW & T-LBBB & T-CRT \\ \hline
		\\[-0.3cm]
		\textbf{OO} & 522 & 490 & 305 & 305 \\ 
		\textbf{OA} & 0 & 27 & 142 & 108 \\ 
		\textbf{A} & 0 & 1 & 23 & 39 \\ 
		\textbf{C} & 9 & 13 & 61 & 79 \\ 
		\hline
	\end{tabular}
	\caption{\ac{PMJs} ($\#$531) classification returned by Algorithm~\ref{alg:pseudocode_bis} in all the test cases. OO=orthodromic activated by AV node; OA=orthodromic activated by antidromic \ac{PMJs}; A=antidromic; C=collision.}
	\label{table:result5}
\end{table}

\subsection{T-H: healthy propagation}
\label{subsec:T-H}
We present a simulation characterized by healthy propagation, which will be used as a baseline for comparisons with the pathological scenarios. Figure~\ref{fig:result1}(b,c) shows the resulting activation map. Only orthodromic propagation is present, following His-Purkinje activation, which takes approximately \SI{30}{\milli\second}, in alignment with the values recorded in-vivo~\cite{verzicco2022}. The sub-endocardial layer is the first portion of the muscle to activate, with the signal progressively spreading throughout the myocardium and transmurally towards the epicardium. The propagation follows an apico-basal pattern, favored by the scarce presence of Purkinje fibers near the ventricular bases. The total ventricular activation times are reported to be around $\SIrange[range-phrase=-]{80}{100}{\milli\second}$ in healthy patients~\cite{palamara_ortho,durrer1970}. The latest activation time, recorded in Table~\ref{table:result5}, is slightly elevated, likely due to the inclusion of the valvular plane in the geometry and to muscle thickening caused by the patient's hypertrophy.

\subsection{T-WPW: Wolff-Parkinson-White syndrome}
The Wolff-Parkinson-White (WPW) syndrome is a pathology characterized by the presence of an accessory pathway between the left atrium and the left ventricle, named bundle of Kent~\cite{sealy1978,al1999clinical,mittal2023}. Since we only consider a ventricular geometry, we follow the modeling approach proposed in~\cite{palamara_anti}, which surrogates the effect of WPW syndrome on the ventricles by introducing an intramuscular source at the termination of the Kent bundle. This induces a pre-excitation of the left ventricle.

The simulation pacing setup is shown in Figure~\ref{fig:result1}(d). A muscular source was placed in the posterolateral portion, near the basal plane, of the left ventricular myocardium, consistently with in-vivo observations of the location of the Kent bundle~\cite{palamara_anti,waller4part}. Here, a spherical stimulus with radius \SI{0.5}{\milli\meter} was applied at time $t=$~\SI{-30}{\milli\second} to replicate the pre-excitation of the muscle, considering that the activation of the AV node starts later at $t=$~\SI{0}{ms}, which is used as reference time across all the scenarios. 

Figure~\ref{fig:result1}(e,f) shows the the resulting activation map.
A front propagating through the muscle is clearly observed before the signal starts to propagate through the network, causing the pre-excitation of the upper posterolateral portion of the left ventricle. At the orthodromic \ac{PMJs}, focal activation is still visible. Away from the muscular source, all the \ac{PMJs} are orthodromic, while near it antidromic or collision \ac{PMJs} arise, as we can see in Figure~\ref{fig:result4}(b). 

\begin{figure}
	\centering
	\includegraphics[width=1\textwidth]{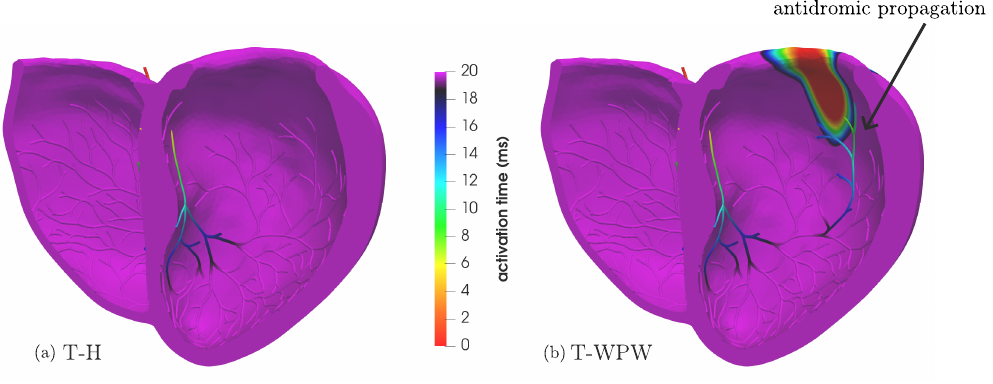}
	\caption{Comparison between T-H (a) and T-WPW (b). The activation time scale was saturated in the range~$\SIrange[range-phrase=-]{0}{20}{\milli\second}$ in order to highlight the identification, in the pathological scenario, of two opposite fronts propagating in the Purkinje network. The effect of antidromic propagation due to the intramuscular source is visible in T-WPW, while only orthodromic propagation is present in T-H.} 
	\label{fig:result2}
\end{figure}

Figure~\ref{fig:result2} reports a comparison between T-H and T-WPW, highlighting the effects of antidromic propagation. We observe that, in the presence of only orthodromic activation, in Figure~\ref{fig:result2}(a), the signal reaches only the three branches diverging from the main bundle of the Purkinje network. Meanwhile, thanks to the pre-excitation of the muscle starting at $t=-$\SI{30}{\milli\second}, the signal enters the network antidromically in Figure~\ref{fig:result2}(b). Here, two wavefronts are distinguishable in opposite directions. Close to the ectopic stimulation point, all the types of \ac{PMJs} are present: antidromic, collision and orthodromic activated by antidromic propagation (OA-labelled \ac{PMJs}). 

\subsection{T-LBBB: complete left bundle branch block}
\label{subsec:T-LBBB}

The second pathological scenario considered is a complete left bundle branch block (LBBB). The latter consists in the interruption of the electric signal propagation in a segment of the left branch, resulting in an alteration of ventricular synchrony~\cite{romano2010,mcanulty1984}. 

The atrioventricular node was stimulated as in the healthy case at time $t=$~\SI{0}{\milli\second} (see Section~\ref{subsec:T-H}). A conduction block is located along a segment of the left bundle branch, shown in black in Figure~\ref{fig:result1}(g).

The total activation time is \SI{174}{\milli\second}, as reported in Table~\ref{table:postpro2} and shown in Figure~\ref{fig:result1}(h,i). The signal originates from the atrioventricular node and travels through the His bundle. At the conduction block, the signal path in the left branch is interrupted, continuing toward the ventricle solely through the right branch. In Figure~\ref{fig:result1}(h), we can notice how the right ventricle undergoes normal orthodromic activation, while the signal in the left ventricle arrives later. The wavefront reaches the left ventricle through the interventricular septum, see Figure~\ref{fig:result1}(h), coherently with what is observed in patients with complete branch blocks. Consequently, impulses enter in the left branch of network antidromically. Figure~\ref{fig:result4}(c) shows that the conduction block determines the presence of antidromic \ac{PMJs} near the interventricular septum. Moreover, no orthodromic \ac{PMJs} activated by the signal coming from the AV node (OO-labelled \ac{PMJs}) are present in the left ventricle. In contrast, we see many OA-labelled \ac{PMJs}, which result from a signal that has entered the network antidromically from the muscle and, due to the higher conduction velocity in the Purkinje fibers, is able to reemerge orthodromically from a different PMJ. A high number of OA \ac{PMJs} is confirmed also by Table~\ref{table:result5}. 

\begin{figure}
	\centering
	\includegraphics[width=1\textwidth]{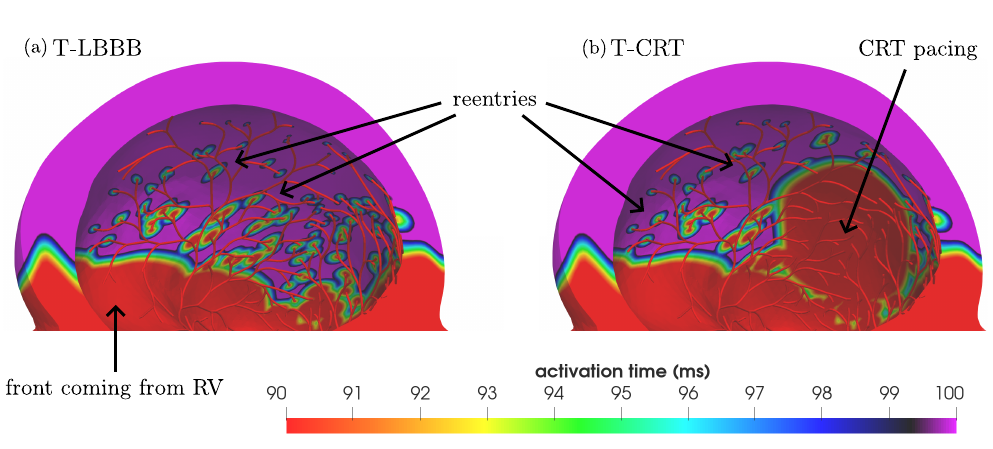}
	\caption{A zoom on the reentries of the electrical signal in the muscle in T-LBBB (a) and T-CRT (b): the muscular wavefront coming from the healthy right ventricle (RV) traverses the interventricular septum and antidromically activates the Purkinje network in the left ventricle, where propagation is interrupted in correspondence of the branch block. Purkinje fibers conduct at higher velocity so that the signal is able to re-emerge orthodromically from \ac{PMJs} (reentries highlighted in picture). To show the fronts, the color scales have been properly saturated.}
	\label{fig:result3}
\end{figure}

Figure~\ref{fig:result3}(a) displays the activation map for the time frame $\SIrange[range-phrase=-]{90}{100}{\milli\second}$, highlighting the reentry propagation from the muscle to the network. The higher transmission velocity in the network causes the signal to re-emerge again in the muscle. In Figure~\ref{fig:result3}(a), a primary front is coming from the right ventricle, overcome by the stimuli arriving from the terminations of the left network. We emphasize that our new proposed Algorithm~\ref{alg:pseudocode_bis} made it possible to numerically observe this reentry effect. In contrast, the methodologies outlined in~\cite{palamara_anti} would not be able to capture such reentries, providing a less biophysically accurate solution. 

\subsection{T-CRT: cardiac resynchronization therapy}
The last presented scenario combines the two sources of antidromic propagation, both the conduction block and the muscular source, into a single simulation, showcasing the versatility of the algorithm for Purkinje-muscle coupling. This test case aims to simulate the effects of \acf{CRT} in presence of a complete left bundle branch block. \ac{CRT} is a treatment that involves biventricular pacing to help restoring the normal rhythm and timing of the heartbeat~\cite{strocchi2022,lee2018,leyva2014}. This involves implanting electrodes in the heart to produce artificial stimuli that compensate for the propagation imbalance caused by a particular pathology. Typically, one lead is placed at the apex of the right ventricle, and another on the free wall midway up the left ventricle~\cite{iaizzo,leyva2014,butter2000,mortensen2004}. The effect of artificial pacing is modeled by defining two appropriate epicardial stimuli~\cite{strocchi2022,lee2018}.

Figure~\ref{fig:result1}(j) resumes the initial setup of the simulation. A stimulus at time $t=$~\SI{0}{\milli\second} is imposed at the AV node, with a complete conduction block located as in T-LBBB (see Section~\ref{subsec:T-LBBB}). Two additional spherical stimuli are imposed in order to mimic the two leads of \ac{CRT}. Activation at these points occurs at time $t=$~\SI{60}{\milli\second} and $t=$~\SI{40}{\milli\second}, consisting of spherical stimuli with radius \SI{0.5}{\milli\meter}. Notice that this is an arbitrary choice, as is their location, which usually corresponds to the termination of two cardiac veins. The correct positioning and timing may vary a lot from patient to patient and are evaluated in clinical practice on individual specifics.

The biventricular activation map is depicted in Figure~\ref{fig:result1}(k,l). The CRT lead in the right ventricle has, on purpose, a higher activation time than the one ensured by orthodromic propagation, in order to verify that it produces no effect. Indeed, the novel pseudo-time method, proposed in Section~\ref{subsec:novel_psuedo}, correctly operates in this situation. Conversely, the lead in the left ventricle produces a stimulus that impacts the solution, as shown in Figure~\ref{fig:result1}(k,l). This is even more visible in Figure~\ref{fig:result3}(b), where the artificial pacing is evident, along with the signal entering the network via antidromic \ac{PMJs}. Several wavefronts are present in the muscle: one main front, which comes from the interventricular septum, and one produced by the \ac{CRT} lead, but also many pointwise stimuli at \ac{PMJs} due to reentries. The right ventricle activates normally, while the left ventricle shows alterations compared to the healthy case. The latter are partially compensated by the artificial pacing, ensuring better activation of the lateral wall of the left ventricle. Figure~\ref{fig:result3} also provides a comparison between T-LBBB (a) and T-CRT (b): the anterolateral portion of the left ventricle is activated thanks to \ac{CRT}. In T-CRT, in addition to the antidromic and collision \ac{PMJs} in correspondence to the interventricular septum, further \ac{PMJs} of these types are present near the \ac{CRT} device, as depicted in Figure~\ref{fig:result4}(d). Moreover, many orthodromic \ac{PMJs} activated due to antidromic propagation (OA type) are visible in the free wall of the left ventricle, consistently with the reentry effect.

	\section{Conclusions}
	\label{sec:conclusions}
	In this paper, we studied the interaction of the electrical excitation in the cardiac muscle and the Purkinje network using Eikonal equations, focusing on modeling electrical reentries. We proposed a novel Eikonal-diffusion pseudo-time method capable of correctly imposing active stimuli in the muscle while discarding inactive ones. Furthermore, we introduced a novel Purkinje-muscle coupling partitioned scheme (Algorithm~\ref{alg:pseudocode_bis}) that effectively manages possible reentries between the two domains by iteratively computing activation times and exchanging information at the interface, represented by the \ac{PMJs}. Building on and extending the work presented in~\cite{palamara_anti}, we addressed its limitations, achieving improved versatility and accuracy.

While this work primarily focused on methodological and modeling aspects, we also performed numerical simulations in realistic pathological and therapeutic scenarios. The results underscore the critical role of the Purkinje network and highlight the importance of accurately capturing electrical reentries to better reproduce the pathophysiological processes. The numerical simulations showcased the accuracy, flexibility and robustness of our algorithm, successfully handling arbitrary combinations of Purkinje network, stimulation setups and conduction blocks.

A promising direction for future work is to rigorously define the continuous formulation of the pseudo-time-discrete problem~\eqref{eq:eikonal3D_pseudo_time}, which is not consistent with~\eqref{eq:eikonal3D} and implicitly introduces a novel notion of solution, whose properties warrant further investigation. Additionally, our numerical experiments relied on Purkinje networks reconstructed through rule-based algorithms, rather than patient-specific. However, thanks to the modularity of the proposed scheme, the latter can be included straightforwardly. Finally, investigating pathological scenarios from an electro-mechanical perspective may provide deeper insights into their effects on cardiac function.
	
	\section*{Acknowledgements}
M.B., R.P. and L.D. received support from the project PRIN2022, MUR, Italy, 2023--2025 202232A8AN ``Computational modeling of the heart: from efficient numerical solvers to cardiac digital twins''. R.P. has also received support from the INdAM GNCS project CUP E53C23001670001 ``Mathematical models and numerical methods for the construction of cardiac digital twins”. C.V. has been partially supported by: i) the European Union-Next Generation EU, Mission 4, Component 1, CUP: D53D23018770001, under the research project MUR PRIN22-PNRR n.P20223KSS2, ``Machine learning for fluid structure interaction in cardiovascular problems: efficient solutions, model reduction, inverse problems'', ii) the Italian Ministry of Health within the PNC PROGETTO HUB LIFE SCIENCE - DIAGNOSTICA AVANZATA (HLS-DA) ``INNOVA'', PNCE3-2022-23683266–CUP: D43C22004930001, within the ``Piano Nazionale Complementare Ecosistema Innovativo della Salute'' - Codice univoco investimento: PNCE3-2022-23683266; iii) the Italian research project MUR PRIN22 n.2022L3JC5T ``Predicting the outcome of endovascular repair for thoracic aortic aneurysms: analysis of fluid dynamic modeling in different anatomical settings and clinical validation''; iv) Italian Ministry of Health within the project ``CAL.HUB.RIA'' - CALABRIA HUB PER RICERCA INNOVATIVA ED AVANZATA. Code: T4-AN-09, CUP: F63C22000530001. M.B., R.P., L.D. and C.V. acknowledge their membership to INdAM group GNCS - Gruppo Nazionale per il Calcolo Scientifico (National Group for Scientific Computing, Italy). M.B. and R.P. acknowledge the INdAM GNCS project CUP E53C23001670001.  M.B., R.P. and L.D. acknowledge the “Dipartimento di Eccellenza 2023-2027”, MUR, Italy, Dipartimento di Matematica, Politecnico di Milano.
\section*{Declaration of competing interest}
All the authors declare that they have no known competing financial interests or personal relationships that could have appeared to influence the work reported in this paper.

\section*{CRediT authorship contribution statement}
\textbf{S. Brunati}: Conceptualization, Methodology, Software, Simulation, Data curation, Formal analysis, Investigation, Visualization, Writing – original draft, Writing – review \& editing.
\textbf{M. Bucelli}: Conceptualization, Methodology, Software, Formal analysis, Investigation, Visualization, Writing – original draft, Writing – review \& editing.
\textbf{R. Piersanti}: Conceptualization, Methodology, Formal analysis, Investigation, Visualization, Writing – original draft, Writing – review \& editing.
\textbf{L. Dede’}: Supervision, Funding acquisition, Project administration, Writing – review \& editing.
\textbf{C. Vergara}: Conceptualization, Supervision, Funding acquisition, Project administration, Writing – original draft, Writing – review \& editing.
	
	\vspace{6mm}
	
	%

\bibliographystyle{elsarticle-num}
\bibliography{references}

\end{document}